\documentclass[10pt]{article}

\newtheorem{proposition}{Proposition}  
\newtheorem{lemma}{Lemma}  
\usepackage{vmargin}   
\setmarginsrb{3.2cm}{2.5cm}{3.2cm}{2.5cm}{1cm}{1cm}{1cm}{1.6cm}   
\usepackage{amsmath}   
\usepackage{amssymb}

\title{An infinite product formula for  $U_q(sl(2))$  dynamical 
coboundary element.}

\author{E.Buffenoir\thanks{e-mail:buffenoi@lpm.univ-montp2.fr}, Ph. Roche\thanks{e-mail:roche@lpm.univ-montp2.fr}{ }\thanks{Work supported by CNRS} \\  
 Laboratoire de Physique Math\'ematique et Th\'eorique  \\  
 Universit\'e Montpellier 2, 34000 Montpellier, France.}  
\date{\today}

\begin{document} 

\maketitle  

\begin{abstract}
We give a short summary of results and conjectures in the theory of dynamical quantum group related to the dynamical coboundary equation also known as IRF-Vertex transform. 
O.Babelon has shown that 
the dynamical twist $F(x)$  of $U_q(sl(2))$ is  a dynamical coboundary $M(x)$ i.e $F(x)M_1(xq^{h_2})M_2(x)=\Delta(M(x)).$
We give a new formula for this element $M(x)$  as an infinite product and give a new proof of the coboundary relation. Our proof involves the quantum Weyl group element, giving possible hint for the generalization to higher rank case.

\end{abstract}

\section*{I. Introduction}
While participating to the conference ``Recent Advances in the Theory of Quantum Integrable Systems 2003'', LAPTH, Annecy-Le-Vieux, we have seen that there was some interest on very preliminary work we have done on the understanding of the shifted coboundary  element some years ago. We have therefore decided to publish these results in the proceedings of the conference.

The theory of dynamical quantum group is nowdays a well established part of mathematics, see the review by P.Etingof \cite{E}. This theory originated from the notion of dynamical Yang-Baxter equation, which arose in the work of Gervais-Neveu on  Liouville theory
\cite{GN} and was formalized  first by G.Felder \cite{Fe} who also understood its relation with IRF statistical models.

In its work on quantum Liouville on a lattice \cite{Bab}, O.Babelon was the first to understand universal aspects of the dynamical Yang-Baxter equation. He  introduced the notion of dynamical twist $F(x)\in U_q(sl(2))^{\otimes 2}$ and  gave  an exact formula for $F(x)$. Quite remarkably he obtained the explicit  formula on $F(x)$ by noting that $F(x)$ is a dynamical coboundary, i.e there exists an explicit invertible element $M(x)\in U_q(sl(2))$ such that:
\begin{equation}
F(x)=\Delta(M(x)) M_2(x)^{-1}(M_1(xq^{h_2}))^{-1}.\label{coboundaryM}
\end{equation}

Twelve years after this work, the theory of dynamical quantum groups is now well understood for $U_q({\mathfrak g})$ where ${\mathfrak g}$ is a Kac-Moody algebra of affine or finite type. 
An important result is the so-called ``linear equation'' \cite{BR}\cite{ABRR}\cite{JKOS} which allows to express  the dynamical twist  as an infinite product of ``dressed'' $R$-matrices. It is quite surprising that the theory of dynamical quantum groups has made no significant progress on the understanding of the coboundary element. See however the articles \cite{BBB}\cite{EN}\cite{St}. 

\medskip 
We will now recall the present situation as well as conjectures. We will give in the next section a new formula for $M(x)$ expressed as  an infinite product. Using this formula we will show that we can recover (after some work !) the infinite product formula for $F(x).$ The present work stems from an unsuccessful attempt to understand the structure of the coboundary for higher rank  case. However we hope that our work  will trigger some further study  of the coboundary element and its  relation with the  quantum Weyl group.

We first recall the results of Babelon using its notations.

In the sequel we will denote $[z]_q=\frac{q^z-q^{-z}}{q-q^{-1}}$ as well as
 $(z)_q=q^{z-1}[z]_q.$
The $q-$exponential is defined as:
\begin{equation}
\exp_q(z)=\sum_{n=0}^{+\infty}\frac{z^n}{(n)_q !}
\end{equation}
where $(n)_q !=(n)_q\cdots (1)_q.$

It satisfies:
\begin{equation}
\exp_q(z)\exp_{q^{-1}}(-z)=1,\label{expinverse}
\end{equation}
and if $xy=q^2yx$
\begin{eqnarray}
&&\exp_q(x+y)=\exp_q(y)\exp_q(x),\label{expsomme}\\
&&\exp_q(x) \exp_q(y)=\exp_q(y)\exp_q((1-q^{-2})xy)\exp_q(y
)\label{expproduit}.
\end{eqnarray}

Let us define $U_q(sl(2))$ to be the Hopf algebra generated by $E_+, E_-, h$ satisfying:
\begin{eqnarray}
&&[h,E_{\pm}]=\pm 2 E_{\pm}, [E_+,E_-]=[h]_q\\
&&\Delta(E_{\pm})=E_{\pm}\otimes q^{h/2}+q^{-h/2}\otimes E_{\pm}.
\end{eqnarray}
We will define $e= q^{h/2}E_+, f=E_- q^{-h/2}.$

We can express the Drinfeld universal $R$ matrix as

\begin{equation}
R=q^{\frac{h\otimes h}{2}}\hat{R}
\end{equation}
 with $\hat{R}=\exp_{q^{-1}}((q-q^{-1})e\otimes f).$

The dynamical twist $F$ is a function $F:{\mathbb C}^{*}\rightarrow U_q(sl(2))^{\otimes 2}$  of zero weight, i.e $[F(x),\Delta(h)]=0$, satisfying the dynamical cocycle equation:
\begin{equation}
F_{12,3}(x)F_{12}(xq^{h_3})=F_{1,23}(x)F_{23}(x),
\end{equation}
where we have denoted $F_{12,3}(x)=(\Delta\otimes id)(F(x)).$

An explicit solution was constructed by O.Babelon and reads:

\begin{equation}
 F(x)=\sum_{k=0}^{+\infty}\frac{(-1)^k (q-q^{-1})^k}{(k)_{q^{-1}}!}
\frac{q^{k(1-k)}}{\prod_{p=1}^k(1-x^{-2}q^{-2p-2h_{2}})}e^k\otimes f^k.
 \label{seriesforF}
\end{equation}

If we define $R(x)=F_{21}(x)^{-1}R F_{12}(x)$, $R(x)$ is a solution of the dynamical Yang-Baxter equation i.e:
\begin{equation}
R_{12}(x)R_{13}(xq^{h_2})R_{23}(x)=R_{23}(xq^{h_1})R_{13}(x)R_{12}(xq^{h_3}).
\end{equation}.

It was shown in \cite{BR}, that $F(x)$ satisfies the linear equation
\begin{equation}
F(x)B_2(x)=\hat{R}^{-1}B_2(x)F(x),
\end{equation}
where $B(x)=x^h q^{h^2/2}.$

This equation can be iterated and 
 we obtain an expression of $F(x)$ in term of an infinite product \cite{ABRR}:
\begin{eqnarray}
&&F(x)=\prod_{k=0}^{+\infty}(B_2(x)^k \hat{R}^{-1} B_2(x)^{-k})\\
&&=\prod_{k=0}^{+\infty} \exp_{q}(-(q-q^{-1})x^{-2k}q^{-2k(1+h_2)}e\otimes f)\label{productforF}
\end{eqnarray}

This product is shown to be convergent in finite dimensional representations and for $x$ sufficiently large.
We can also express $F(x)$ as:
\begin{eqnarray}
F(x)=\prod_{k=0}^{+\infty}(B_2(x)^{-k-1} \hat{R} B_2(x)^{k+1}),
\end{eqnarray}
which is convergent  in finite dimensional representations and for $x$ sufficiently small.

Babelon has shown \cite{Bab} that if we define $M:{\mathbb C}^{*}\rightarrow U_q(sl(2))$ by:
\begin{equation}
M(x)=\sum_{n,m=0}^{+\infty }
\frac{(-1)^m x^m q^{n(n-1)/2+m(n-m)}}
{[m]_q! [n]_q!
\prod_{p=1}^{n}(xq^p-x^{-1}q^{-p})
}E_+^n E_-^m q^{(n+m)h/2},\label{Mofx}
\end{equation}

we have 
\begin{equation}
F(x)M_1(xq^{h_2})M_2(x)=\Delta(M(x)).
\end{equation}

Note that  
 $M(x)$ is not of  zero weight, i.e $[M(x),h]\not=0.$

{}From this relation we obtain that:

\begin{equation}
R(x)M_1(xq^{h_2})M_2(x)= M_2(xq^{h_1}) M_1(x)R,
\end{equation} 

i.e one can absorb all the dynamics through the use of $M(x).$

Note that $M(x)$ relates the dynamical $R$-matrix $R(x),$ which in the affine case, would  give the IRF Boltzmann weights, to the $R$ matrix, which in the affine case, would  give the vertex Boltzmann weights. This is the reason why $M(x)$ is also called by extension IRF-Vertex transform. Such type of $M,$ in matrix form, has been computed by Baxter in \cite{Bax} and   relates the IRF-model to the 8-Vertex solution.

All these results have been obtained in the $U_q(sl(2))$ case. Since we now have a good theory of dynamical quantum groups in the $U_q({\mathfrak g})$ case, one can expect that such results can be generalized to higher rank case.

We will use the notations of \cite{ABRR}. We denote $F(x)$ the dynamical twist in the $U_q(\mathfrak{g})$ case where $x=q^{\lambda}$ with $\lambda\in \mathfrak{ h}.$
The linear equation still applies in this case and we have:
\begin{equation}
F(x)=\prod_{k=0}^{+\infty}(B_2(x)^k \hat{R}^{-1} B_2(x)^{-k}),\label{2shiftedcocycle}
\end{equation}
where
$B(x)$ and  $\hat{R}$ is the straightforward generalization to higher rank as explained in \cite{ABRR}.

The product is convergent in finite dimensional representations and for $\lambda$ sufficiently large.

The first result on coboundary element for higher rank case was obtained  by Cremmer-Gervais \cite{CG} in the $sl(3)$ case and by  Bilal-Gervais \cite{BG} in the $sl(n+1)$ case, where they have shown that, in the fundamental representation of  $sl(n+1),$ one can absorb all the dynamic of $R(x)$ and one obtains:

\begin{equation}
R(x)M_1(xq^{h_2})M_2(x)= M_2(xq^{h_1}) M_1(x)R_{CG},
\end{equation} 
with $M(x)$ an $(n+1)\times (n+1)$ explicit matrix (of VanderMonde type)  and 
$R_{CG}$ satisfying Yang-Baxter equation.

Quite surprisingly $R_{CG}$  is not the standard solution of Yang-Baxter equation and is not of zero weight.
Its classical expansion $R_{CG}=1+\hbar r_{CG}+o(\hbar)$ is associated to a particular solution of Belavin-Drinfeld classification.
In the $sl(n+1)$ case it corresponds to taking the following Belavin-Drinfeld triple:
\begin{eqnarray*}
&&\Gamma=\{\alpha_1,\cdots,\alpha_{n}\}\\
&&\Gamma_1=\{\alpha_1,\cdots,\alpha_{n-1}\},\Gamma_2=
\{\alpha_2,\cdots,\alpha_{n}\}\\
&&T:\Gamma_1\rightarrow \Gamma_2, \alpha_i\mapsto\alpha_{i+1}.
\end{eqnarray*}

Unfortunately all these computations have been done in the fundamental representations and gives no real hints towards a  precise universal formulation.
The subject has seen a dramatic advance with the result of Etingof, Schedler, Schiffmann \cite{ESS}, where they have shown how to quantize explicitly all classical  solutions of Yang-Baxter equation given by a Belavin-Drinfeld triple.
They have constructed for all Belavin-Drinfeld triple $T$ a twist $J_T=J\in U_q({\mathfrak g})^{\otimes 2}$, which satisfies 

\begin{equation}
J_{12,3}J_{12}=J_{1,23}J_{23}.
\end{equation}

and is such that $J_{21}^{-1}RJ_{12}$ satisfies Yang-Baxter equation and is a
 quantization of the classical $r_T$ matrix associated to $T.$
The expression of $J$ was obtained through a nice use of dynamical quantum groups and of a modification of the linear equation. Finally  $J$ is expressed as a finite product of explicit invertible elements.
Therefore this result provides the answer of the construction of explicit universal formulas for the non standard solution of Yang-Baxter equation such as the Cremmer-Gervais 's one. 

It is therefore tempting to formulate the following conjectures:
\medskip

{\bf Conjecture 1} When ${\mathfrak g}=sl(n+1),$
there exists $M:{\mathbb C}^{*n}\rightarrow U_q({\mathfrak g})$ such that:
\begin{equation}
F(x)M_1(xq^{h_2})M_2(x)=\Delta(M(x))J_{T},\label{coboundarysln}
\end{equation}
with $T$ given by the previous Belavin-Drinfeld triple. 

\medskip

{\bf Conjecture 2 }A solution to this equation is given in term of an infinite product of explicit invertible elements of $U_q({\mathfrak g})$

In the sequel we will reconsider these conjectures in the $U_q(sl(2))$ case.
 The conjecture 1, in this case is already solved, but we will give a new formula for $M(x)$ in this case and show that conjecture 2 is satisfied. We will then  check directly on the infinite product expression that  the coboundary 
equation is satisfied.

\section*{II.The infinite product expression.}

Let $Q:{\mathbb C}^{*}\rightarrow U_q(sl\; 2).$
We denote $\delta Q$ the function 
$\delta Q(x)=\Delta(Q(x))Q_2(x)^{-1}(Q_1(xq^{h_2}))^{-1}.$

A dynamical group like element is a function $g:{\mathbb C}^{*}\rightarrow U_q(sl\; 2)$ such that $\delta g(x)=1.$

One can easily construct such dynamical group like elements.
Apart the constant group like element, one can also define for example: 
\begin{eqnarray}
\text{Example 1}&&g(x)=B(x)\\
\text{Example 2}&&g(x)=(x;q)_h=\frac{(x;q)_{\infty}}{(xq^h;q)_{\infty}}.
\end{eqnarray}

If $g:{\mathbb C}^{*}\rightarrow  U_q(h),$ is a dynamical group like element, 
 then $$\delta (gQ)(x)=\Delta(g(x))\delta(Q(x))(\Delta(g(x)))^{-1}.$$
As a result if $Q$ satisfies the coboundary equation (\ref{coboundaryM}), from the $h$ invariance of $F(x)$, we obtain that $gQ$ is also a solution of the coboundary equation.

We will now show that there exists a solution of the coboundary equation which is written as a simple infinite product formula.  We think that this proof, although a bit intricate, could be generalized to higher rank case in spite of unsuccessful attempt from us.

\begin{proposition}
Let $N:{\mathbb C}^*\rightarrow U_q(sl 2)$ defined as:
\begin{equation}
N(x)=N_-(x)N_+(x)
\end{equation}
with
\begin{eqnarray}
&&N_+(x)=\exp_{q^{-1}}(-x e)\\
&&N_-(x)=\prod_{k=+\infty}^0\exp_{q^{-1}}(q^{-(2k+1)(h+1)}x^{-2k-1}f).
\end{eqnarray}
$N(x)$  satisfies the coboundary equation.
\end{proposition}

In order to prove this proposition, we first begin by proving a lemma interesting in itself.

We use the notations of the article \cite{ABRR}: 
let $\mathfrak{g}$ be a finite simple Lie algebra of rank $r$ and denote by 
$F:{\mathbb C}^{*r}\rightarrow U_q(\mathfrak{g})^{\otimes 2}$  the solution of the dynamical cocycle equation.

\begin{lemma}
Let $u,v$ be maps from 
${\mathbb C}^{*r}$ to $U_q(\mathfrak{g})$ and ${\cal J}$ an invertible element of $U_q(\mathfrak{g})^{\otimes 2}$.

We define:
\begin{eqnarray*}
V^{(p)}(x)=\prod _{k=p}^{1}\left(B(x)^k v(x)B(x)^{-k}\right),V^{(0)}(x)=1  &  & M^{(0)}(x)= u(x)^{-1} \\
M^{(p)}(x)=V^{(p)}(x)M^{(0)}(x) &  & F^{(p)}(x)= \Delta \left(M^{(p)}(x)\right)\: \mathcal{J}\: M^{(p)}(x)_{2}^{-1}M^{(p)}(xq^{h_{2}})_{1}^{-1}.
\end{eqnarray*}
we define $G(x)\equiv F^{(0)}(x)^{-1}\Delta(B(x))^{-1}F^{(1)}(x)\Delta(B(x))$.

If the following relation is satisfied:
\begin{equation}
\left[G(x)\: v(xq^{h_{2}})_{1}B(xq^{h_{2}})_{1}^{-1}\: ,\: M^{(0)}(xq^{h_{2}})_{1}v(x)_{2}B(x)_{2}^{-1}M^{(0)}(xq^{h_{2}})_{1}^{-1}\right]=0
\label{commutationG1}
\end{equation}
as well as the asymptotic relations:
\begin{eqnarray*}
\lim _{k\rightarrow +\infty }B(x)^{k}v(x)^{-1}B(x)^{-k} & = & 1\\
\lim _{k\rightarrow +\infty }\Delta \left(B(x)\right)^{k}F^{(0)}(x)\Delta \left(B(x)\right)^{-k} & = & 1\\
\lim _{k\rightarrow +\infty }\Delta \left(B(x)\right)^{k}G(x)\Delta \left(B(x)\right)^{-k} & = & \hat{R}^{-1},
\end{eqnarray*}
then 
\begin{eqnarray*}
\lim _{k\rightarrow +\infty }F^{(k)}(x) & = & F(x).
\end{eqnarray*}
\end{lemma}
Proof:
We first remark the following recursion relations:

\begin{eqnarray*}
&&V^{(p)}(x)=B(x)V^{(p-1)}(x)v(x)B(x)^{-1}\\
&&  M^{(p)}(x)=B(x)M^{(p-1)}(x)\left(M^{(0)}(x)^{-1}v(x)B(x)^{-1}M^{(0)}(x)\right).
\end{eqnarray*}
This  allows us to deduce the following recursion relation on $F^{(p)}(x):$
\begin{eqnarray*}
F^{(p)}(x) & = & \Delta \left(B(x)\right)\Delta \left(M^{(p-1)}(x)\right)\Delta \left(M^{(0)}(x)^{-1}v(x)B(x)^{-1}M^{(0)}(x)\right)\: \mathcal{J}\: M^{(p)}(x)_{2}^{-1}M^{(p)}(xq^{h_{2}})_{1}^{-1}\\
 & = & \Delta \left(B(x)\right)F^{(p-1)}(x)\left(M^{(p-1)}(xq^{h_{2}})_{1}M^{(p-1)}(x)_{2}\: \mathcal{J}^{-1}\: \Delta \left(M^{(0)}(x)^{-1}v(x)B(x)^{-1}M^{(0)}(x)\right)\times \right.\\
 &  & \left.\times \: \mathcal{J}\: M^{(0)}(x)_{2}^{-1}V^{(p)}(x)_{2}^{-1}M^{(0)}(xq^{h_{2}})_{1}^{-1}V^{(p)}(xq^{h_{2}})_{1}^{-1}\right)\\
 & = & \Delta \left(B(x)\right)F^{(p-1)}(x)\left(M^{(p-1)}(xq^{h_{2}})_{1}M^{(p-1)}(x)_{2}M^{(0)}(x)_{2}^{-1}M^{(0)}(xq^{h_{2}})_{1}^{-1}\: G(x)\: v(xq^{h_{2}})_{1}B(xq^{h_{2}})_{1}^{-1}\times \right.\\
 &  & \left.\times M^{(0)}(xq^{h_{2}})_{1}V^{(p-1)}(x)_{2}^{-1}B(x)_{2}^{-1}M^{(0)}(xq^{h_{2}})_{1}^{-1}B(xq^{h_{2}})_{1}v(xq^{h_{2}})_{1}^{-1}V^{(p-1)}(xq^{h_{2}})_{1}^{-1}B(xq^{h_{2}})_{1}^{-1}\right).
\end{eqnarray*}
In order to check the last equality, it is sufficient to show the following equality after  having expanded $V^{(p)}(x)$ in terms of  $V^{(p-1)}(x):$

\begin{eqnarray*}
&&\mathcal{J}^{-1}
\Delta (M^{(0)}(x)^{-1}v(x)B(x)^{-1}M^{(0)}(x))
\mathcal{J} M^{(0)}(x)_{2}^{-1}V^{(p)}(x)_{2}^{-1}=\\
&& M^{(0)}(x)_{2}^{-1}M^{(0)}(xq^{h_{2}})_{1}^{-1}\: G(x)\: v(xq^{h_{2}})_{1}
B(xq^{h_{2}})_{1}^{-1} M^{(0)}(xq^{h_{2}})_{1}V^{(p-1)}(x)_{2}^{-1}
B(x)_{2}^{-1}.
\end{eqnarray*}

The left handside of this equation is equal to:

\begin{eqnarray*}
&& M^{(0)}(x)_{2}^{-1}M^{(0)}(xq^{h_{2}})_{1}^{-1}
F^{(0)}(x)^{-1}\Delta(B(x)^{-1})F^{(1)}(x)\Delta(B(x))\times\\
&&\Delta(B(x)^{-1})M^{(1)}(xq^{h_2})_1M^{(1)}(x)_2M^{(0)}(x)_2^{-1}V^{(p)}(x)_2^{-1}\end{eqnarray*}

which is shown to be equal to the righthandside of this equation using the definition of $G(x)$ and the equality,
\begin{eqnarray*}
&&\Delta(B(x)^{-1})M^{(1)}(xq^{h_2})_1M^{(1)}(x)_2M^{(0)}(x)_2^{-1}V^{(p)}(x)_2^{-1}=\\
&& v(xq^{h_{2}})_{1}
B(xq^{h_{2}})_{1}^{-1} M^{(0)}(xq^{h_{2}})_{1}V^{(p-1)}(x)_{2}^{-1}
B(x)_{2}^{-1}.
\end{eqnarray*}

If $G(x)$ satisfies the  commutation relation of the hypothesis it necessarily satisfies: \begin{equation}
\left[G(x)\: v(xq^{h_{2}})_{1}B(xq^{h_{2}})_{1}^{-1}\: ,\: M^{(0)}(xq^{h_{2}})_{1}B(x)_{2}v(x)_{2}^{-1}M^{(0)}(xq^{h_{2}})_{1}^{-1}\right]=0\label{commutationG}\end{equation}
then, by recursion, we deduce 
\[
\left[G(x)\: v(xq^{h_{2}})_{1}B(xq^{h_{2}})_{1}^{-1}\: ,\: M^{(0)}(xq^{h_{2}})_{1}V^{(p-1)}(x)_{2}^{-1}B(x)_{2}^{p-1}M^{(0)}(xq^{h_{2}})_{1}^{-1}\right]=0.\]

We will now show that this relation implies:

\begin{eqnarray*}
&&(M^{(p-1)}(xq^{h_{2}})_{1}M^{(p-1)}(x)_{2}M^{(0)}(x)_{2}^{-1}M^{(0)}(xq^{h_{2}})_{1}^{-1}\: G(x)\: v(xq^{h_{2}})_{1}B(xq^{h_{2}})_{1}^{-1}\times \\
 &  & M^{(0)}(xq^{h_{2}})_{1}V^{(p-1)}(x)_{2}^{-1}B(x)_{2}^{-1}M^{(0)}(xq^{h_{2}})_{1}^{-1}B(xq^{h_{2}})_{1}v(xq^{h_{2}})_{1}^{-1}V^{(p-1)}(xq^{h_{2}})_{1}^{-1}B(xq^{h_{2}})_{1}^{-1}=\\
&&V^{(p-1)}(xq^{h_{2}})_{1}B(x)_{2}^{p-1} G(x) (V^{(p-1)}(xq^{h_{2}})_{1}B(x)_{2}^{p-1})^{-1}\Delta \left(B(x)\right)^{-1}.\\
\end{eqnarray*}

Indeed in this last equation by eliminating $V^{(p-1)}(xq^{h_{2}})_{1}$ on the left and eliminating $(V^{(p-1)}(xq^{h_{2}})_{1}B(x)_{2}^{p-1})^{-1}$ on the right, we obtain the previous relation.

It is now easy to simplify the recursion relation for $F^{(p)}(x):$\begin{eqnarray*}
F^{(p)}(x) & = & \Delta \left(B(x)\right)F^{(p-1)}(x)\left(V^{(p-1)}(xq^{h_{2}})_{1}B(x)_{2}^{p-1}\right)G(x)\left(V^{(p-1)}(xq^{h_{2}})_{1}B(x)_{2}^{p-1}\right)^{-1}\Delta \left(B(x)\right)^{-1}.
\end{eqnarray*}

Using the previous relation, as well as the formula \begin{eqnarray*}
B(x)_{2}^{-i}V^{(i-1)}(xq^{h_{2}})_{1}^{-1}\Delta \left(B(x)\right)^{-1}V^{(i)}(xq^{h_{2}})_{1}B(x)_{2}^{i+1} & = & v(xq^{h_{2}})_{1}B(xq^{h_{2}})_{1}^{-1}
\end{eqnarray*}

we obtain finally \begin{eqnarray*}
F^{(p)}(x) & = & \left(\Delta \left(B(x)\right)^{p}F^{(0)}(x)\Delta \left(B(x)\right)^{-p}\right)\times \\
 &  & \times \prod_{i=0}^{p-1}\left(B(x)_{2}^{i}\left(\Delta \left(B(x)\right)^{p-i}G(x)\Delta \left(B(x)\right)^{-(p-i)}\right)\left(B(xq^{h_{2}})_{1}^{p-i}v(xq^{h_{2}})_{1}B(xq^{h_{2}})_{1}^{-(p-i)}\right)B(x)_{2}^{-i}\right)\times \\
 &  & \times \prod_{i=1}^{p}\left(B(xq^{h_{2}})_{1}^{i}v(xq^{h_{2}})_{1}^{-1}B(xq^{h_{2}})_{1}^{-i}\right).
\end{eqnarray*}

Using the asymptotic properties, we have:
\begin{eqnarray*}
\!\!\!\!&&\lim_{p\rightarrow +\infty }\left(\Delta \left(B(x)\right)^{p}F^{(0)}(x)\Delta \left(B(x)\right)^{-p}\right)=1\\
\!\!\!\!&&\lim_{p\rightarrow +\infty } \prod_{i=0}^{p/2-1}\left(B(x)_{2}^{i}\left(\Delta \left(B(x)\right)^{p-i}G(x)\Delta \left(B(x)\right)^{-(p-i)}\right)\left(B(xq^{h_{2}})_{1}^{p-i}v(xq^{h_{2}})_{1}B(xq^{h_{2}})_{1}^{-(p-i)}\right)B(x)_{2}^{-i}\right)\\
&&=F(x)\\
\!\!\!\!&&\lim_{p\rightarrow +\infty } \prod_{i=p/2}^{p-1}\left(B(x)_{2}^{i}\left(\Delta \left(B(x)\right)^{p-i}G(x)\Delta \left(B(x)\right)^{-(p-i)}\right)\left(B(xq^{h_{2}})_{1}^{p-i}v(xq^{h_{2}})_{1}B(xq^{h_{2}})_{1}^{-(p-i)}\right)B(x)_{2}^{-i}\right)\times\\
\!\!\!\!&&\times \prod_{i=1}^{p/2}\left(B(xq^{h_{2}})_{1}^{i}v(xq^{h_{2}})_{1}^{-1}B(xq^{h_{2}})_{1}^{-i}\right)=1\\
\!\!\!\!&&\lim_{p\rightarrow +\infty }\prod_{i=p/2+1}^{p}\left(B(xq^{h_{2}})_{1}^{i}v(xq^{h_{2}})_{1}^{-1}B(xq^{h_{2}})_{1}^{-i}\right)=1,
\end{eqnarray*}
we obtain the announced result, i.e.
\begin{eqnarray*}
\lim_{p\rightarrow +\infty }F^{(p)}(x) & = & F(x).
\end{eqnarray*}

We now prove the proposition using this lemma.

\medskip
Proof:

We restrict to $U_q(sl(2))$ and we define the following elements:
\begin{equation}
u=\exp_{q}(e), v=\exp_{q^{-1}}(f),  b(x)=q^{\frac{h^{2}}{4}}
x^{\frac{h}{2}}.
\end{equation}

We set $v(x)=b(x)^{-1}vb(x), u(x)=x^{h/2}u x^{-h/2},
$ and set $\mathcal{J}=1.$
They satisfy the following relations:
\begin{eqnarray}
&&\Delta (u)=u_{2\: }q^{\frac{h_{1}h_{2}}{2}}u_{1\: }q^{-\frac{h_{1}h_{2}}{2}}\label{Deltau}\\
&&\Delta (v)=q^{\frac{h_{1}h_{2}}{2}}v_{2}\: q^{-\frac{h_{1}h_{2}}{2}}v_{1}\label{Deltav}\\
&&u_{1}\: q^{\frac{h^{2}_{1}}{4}}R_{12}^{-1}q^{-\frac{h^{2}_{1}}{4}}\: v_{2}^{-1}=q^{-\frac{h_{1}h_{2}}{2}}\: v_{2}^{-1}\: q^{\frac{h_{1}h_{2}}{2}}\: u_{1}\: q^{-\frac{h_{1}h_{2}}{2}}\label{pentagonale}
\end{eqnarray}

Using these relations we verify that 
$F^{(0)}(x)=1, F^{(1)}(x)={\hat R}^{-1}=G(x).$ 

The relation (\ref{commutationG}) is shown to be satisfied using 
(\ref{Deltav}) and (\ref{pentagonale}). 

It remains to show the asymptotic properties. The first one is  satisfied in finite dimensional representations and for $x$ sufficiently large, the second and the third ones are trivial because the sequences are constant.
$\Box.$

\medskip

We have designed  the proof of the coboundary equation having in mind to generalize it to higher rank case.
The important relation (\ref{commutationG}) was a direct consequence of the 
two relations (\ref{Deltav},\ref{pentagonale}).
We can rewrite these relations using the quantum Weyl group element.
Indeed let $w$ be defined by:
\begin{equation}
w=v q^{-\frac{h^2}{4}} u^{-1}q^{-\frac{h^2}{4}}v= 
u^{-1}q^{-\frac{h^2}{4}}v q^{-\frac{h^2}{4}} u^{-1}.
\end{equation}
This element satisfies 

\begin{eqnarray}
&&whw^{-1}=-h, wew^{-1}=-q^{-h-1}f,  wfw^{-1}=-eq^{h+1}\label{wonx}\\
&&\Delta(w)={\hat R}^{-1}w_1w_2.\label{Deltaw}
\end{eqnarray}

We remark that the pair $(u,v)$ solution of (\ref{Deltau},\ref{Deltav},\ref{pentagonale}) can be constructed solely from $u$ and $w.$
Precisely: if $u$ and $w$ satisfy the relations 
(\ref{Deltau},\ref{wonx}, \ref{Deltaw}) as well as the relation
\begin{equation}
wuw=u^{-1}q^{-\frac{h^2}{2}}wu^{-1},\label{uw}
\end{equation}
then we can define $v=q^{-\frac{h^2}{4}}wu^{-1}w^{-1}q^{\frac{h^2}{4}}$
and the couple $(u,v)$ satisfies the relations  (\ref{Deltau},\ref{Deltav},\ref{pentagonale}).

This is why we strongly suspect that the structure of the proof will remains unchanged in the $U_q(sl(n+1)):$
$w$ will be the longest element in the quantum Weyl group and the only unknown element to us is the element $u$, which definition must  generalize the relations (\ref{Deltau},  \ref{uw}) by incorporating $J_T.$ 

We now come back to the relation between this infinite production solution and the solution of Babelon: they are related by a dynamical group like element.
 It amounts to reorder the $e$ and $f$ in Babelon's formula and to transform a sum in a product. We found it in a very indirect way using the relation between the matrix element of $M(x)$ and the $3j$-symbols of 
$U_q(sl(2))$ \cite{BBB}. 
It was for us very surprising that such a factorization occurs because the formula for $M(x)$ involves the factor $q^{mn}$ which prevents to disentangle the double sum in the opposite factorization. They are related by a dynamical group like element:

\begin{proposition}
Let $M(x)$ be the element defined by (\ref{Mofx}) it satisfies:
$N(x)=(q^2x^2;q^2)_h M(x).$
\end{proposition}

Proof:
We first reorder the $e$ and the $f$ in Babelon's formula.
We will first show that
\begin{eqnarray}
M(x)=\tilde{N}_-(x)N_+(x)
\end{eqnarray}
where 
\begin{equation}
\tilde{N}_-(x)=\sum_{k=0}^{+\infty}
\frac{q^{\frac{k(k+1)}{2}}(-z)^k q^{kh}}{(q^2z^2;q^2)_{h+k}[k]!}f^k.
\end{equation}
By projecting on each weight space for the adjoint action, this amounts to show that the following identities are satisfied for every $r\in \mathbb{Z}:$

\begin{eqnarray}
&&\sum_{n,m\geq 0,n-m=r}\frac{(-x)^mq^{-n+m(m+1)/2}q^{mh}}{[m]![n]!\prod_{p=1}^n(xq^p-x^{-1}q^{-p})}e^nf^m=\label{leftweightr}\\
&&=\sum_{k,l\geq 0,l-k=r}
\frac{q^{\frac{k(k+1)}{2}}(-x)^k q^{kh}}{(q^2x^2;q^2)_{h+k}[k]!}\frac{(-x)^l}{(l)_{q^{-1}}!}f^k e^l.\label{rightweightr}
\end{eqnarray}

To prove this identity it is sufficient to show that it holds in any finite dimensional representation of $U_q(sl(2)).$
The simple $p+1$ dimensional module of $U_q(sl(2)))$ has a basis $v_k$ with $k=0,...,p,$ with the following action:
\begin{equation}
h v_k=(p-2k)v_k,
e v_k=[p-k+1]_qv_{k-1},
 f v_{k}=[k+1]_qv_{k+1}.
\end{equation}

We will show it for $r \geq 0$, the proof is similar for $r<0.$
A straightforward computation shows that the left-handside of (\ref{leftweightr}) acting on $v_k$ is equal to:
\begin{equation}
\frac{[p-k+r]!q^{-r}}{[p-k]![r]!\prod_{p=1}^r(xq^p-x^{-1}q^{-p})}
{}_2\varphi_2(q^{2(k+1)},q^{2(k-p)};q^{2(r+1)}, x^2q^{2(r+1)};q^2)(x^2q^{2(p-2k+2r+1)})v_{k-r}
\end{equation}

whereas the action of the right-handside of (\ref{rightweightr}) acting on $v_k$
gives:
\begin{equation}
\frac{[p-k+r]!(-x)^rq^{\frac{r(r-1)}{2}}}{[p-k]![r]!(q^2x^2;q^2)_{p-2k+2r}}
{}_2\varphi_2(q^{2(p+r-k+1)},q^{2(r-k)};q^{2(r+1)}, x^2q^{2(p-2k+2r+1)};q^2)(x^2q^{2(r+1)})v_{k-r}.\end{equation}

These two expressions are equal	 thanks to the following identity on hypergeometric functions

\begin{equation}
{}_2\varphi_2(a,b;c,d;q)(cd/ab)=\frac{(cd/ab;q)_{\infty}}{(d;q)_{\infty}}
{}_2\varphi_2(c/a,c/b;c,cd/ab;q)(d),
\end{equation}
which is proved using the transformation of ${}_3\varphi_2$ series (formula III.9) of the appendix of the Gasper-Rahman book \cite{GR}.

It remains to show that 
${\tilde N}_{-}(x)=N_-(x),$
this is easily proven by equating the  projection  on the second space   of  the two expressions of $F(x)$ as a series (\ref{seriesforF}) and as a product 
(\ref{productforF}).

\section*{III.Conclusion}

We hope that this result will trigger some further work on this subject.
Thirteen years after the work of Cremmer-Gervais, there is still no universal formula for $M(x)$ for $U_q({\mathfrak g})$ with
 $rank({\mathfrak g})>1.$ This problem has certainly a reasonable answer. From the study of the $U_q(sl(2)$ case it shows that the quantum Weyl group will certainly play a decisive role. 
There are  questions that arose  in the discussions with experts on the subject and which should be answered:

-What is the validity of the conjecture 1?  What is the reason of the appearance of a specific Belavin-Drinfeld  triple? 

-Can it be generalized to other Kac-Moody algebras of finite type or affine type? 

-In the affirmative case, can it be used to define an IRF-Vertex transform in the affine case with potential applications to integrable systems?

-What is the relation, if any, between $M(x)$ and the dynamical quantum Weyl group \cite{EV}? 

\medskip

{\bf Added comment.}
The referee of this article pointed out that H. Rosengren has proved a factorization result of $M(x)$ in \cite{Ro} in a different context.

\medskip

{\bf  Acknowledgments:}

We thank  M.Jimbo, O.Schiffman and A.Varchenko for discussions on the dynamical coboundary element three years ago.   

We warmly thank P.Etingof  for numerous discussions and   collaborating  on this subject three years ago.

We warmly thank O.Babelon   for numerous discussions and convincing  us to publish these results.

We thank the organizers of the conference, D.Arnaudon, J.Avan, L.Frappat, E.Ragoucy, P.Sorba, for the invitation to this conference as well as for the opportunity to publish these results in the proceedings.

We thank the referee for helpful remarks.

\end{document}